\documentclass[english,twoside]{iopart}
\usepackage{epsfig}
\usepackage{graphicx}
\usepackage{amsfonts}
\usepackage{latexsym}
\usepackage{amssymb}
\usepackage{graphicx}

\begin{document}

\title[]{Doppler spectrum estimation by Ramanujan-Fourier Transform (RFT)
}

\author{Mohand Lagha,$^{\dag}$ Messaoud Bensebti$^{\ddag}$ }

\address{$^{\dag}$ Aeronautics Department,
Saad Dahlab University of Blida, B.P. 270,
Road of Soumaa, Blida, Algeria
}

\address{$^{\ddag}$ Electronics Department,
Saad Dahlab University of Blida, B.P. 270,
Road of Soumaa, Blida, Algeria
}

\begin{abstract}

The Doppler spectrum estimation of a weather radar signal in a
classic way can be made by two methods, temporal one based in
the autocorrelation of the successful signals, whereas the other
one uses the estimation of the power spectral density PSD by using
Fourier transforms. We introduces a new tool of signal processing
based on Ramanujan sums $cq(n)$, adapted to the
analysis of arithmetical sequences with several resonances $p/q$.
These sums are almost periodic according to time $n$ of resonances
and aperiodic according to the order $q$ of resonances.
New results will be supplied by the use of Ramanujan Fourier Transform
(RFT) for the estimation of the Doppler spectrum for the weather radar signal.

\end{abstract}


\section{Introduction}

In the signal processing field, transforms or methods have been used to move from a space into another one (time towards frequency) or whatever, in order to estimate and analyze, in a better way, the informational content of the signal. Up to now, Discrete Fourier Transform (DFT) and its Fast Fourier Transform (FFT) have been the best tools ever used for periodic or quasi-periodic signals. But this technique is not really appropriate for the analysis of aperiodic random signals. This fact is not new. As a result, multiple methods have been developed recently to analyze time series, the wavelet method, AR model or the autoregressive moving average (ARMA) model.

In the context of the spectrum parameter estimation of a weather signal radar, an algorithm called Pulse Pair (PP) and based on the calculation of the autocorrelation of complex time series, \textit{Z }(\textit{I,Q), }has been developed at first, and this with one or two lags development of the spectrum moments in Mc-Lauren series in order to reduce the calculation time. Then, thanks to computers and data processors improvement, the Fourier Transform has been introduced, followed by AR model and ARMA to be used in the spectrum domain and analyze, in a better way, the detected weather phenomenons.

This article aims mainly at introducing and studying the Ramanujan sums, \textit{cq}(\textit{n),} in order to estimate the spectrum properties of pulsed Doppler radar signals, used in meteorology.

The use of the Ramanujan method in the domain of the Doppler spectrum estimation of the weather radar  signals, has been motivated by the recent researchers' growing interest to introduce it as a new tool in signal processing [9][10]. This method has been used by Ramanujan, for the first time, as a means for representing arithmetic series by infinite extent sums.  Moreover, transform coefficients,\textit{ cq,} come from the use of arithmetic functions of the number theory.  We may add that these coefficients have the orthogonality property which makes them much more interesting in signal processing.

In this article, the spectrum moment estimation - zero, one and two orders (power, velocity and spectrum width) for a pulsed Doppler weather radar signal (section 2) -- has been dealt with and the two estimation domains: The first one which is based on the autocorrelation estimation (Time domain estimation), and the second one which is based on power spectrum densities PSD (frequency domain estimation), have been treated while using the Discrete Fourier Transform (DFT). A certain notion of the Ramanujan sums has been introduced and the Ramanujan-Fourier Transform (RFT) presented as a new tool in signal processing of arithmetic sequences, based on the prime number theory (section 3). Results about the calculation made on the real data provided by a WSR-88 D pulse Doppler radar and discussions about estimations led upon the Doppler spectrum by the estimators cited (section 4).

\section{Spectral moment estimation}

The radar data, repeated at regular time intervals is referred to as a volume scan. The reflectivity measured in \textit{dBz} often contained in weather Doppler radar signal may contain precipitation and wind information.

The pulsed Doppler radar delivers the output voltages - I in-phase and Q in quadrature phase -- thus, generating the complex echo \textit{Z (I, Q)}. The power spectrum density is given by the Fourier Transform of the autocorrelation function  $R_{ZZ} (\tau )$. [1], [2].

\begin{equation} \
S_{Z} (f)=\Im \left\{R_{ZZ} (\tau )\right\}
\end{equation}

Where \textit{Z (I, Q) =I+jQ }is a complex signal generated at radar receiver by the weather echoes returned [1], from the weather perturbations\textit{. }The received Doppler spectrum $S_{Z} (f)$is represented in figure 1.

\begin{figure}[ht]
\centerline{\includegraphics[width=15.0truecm,clip=]{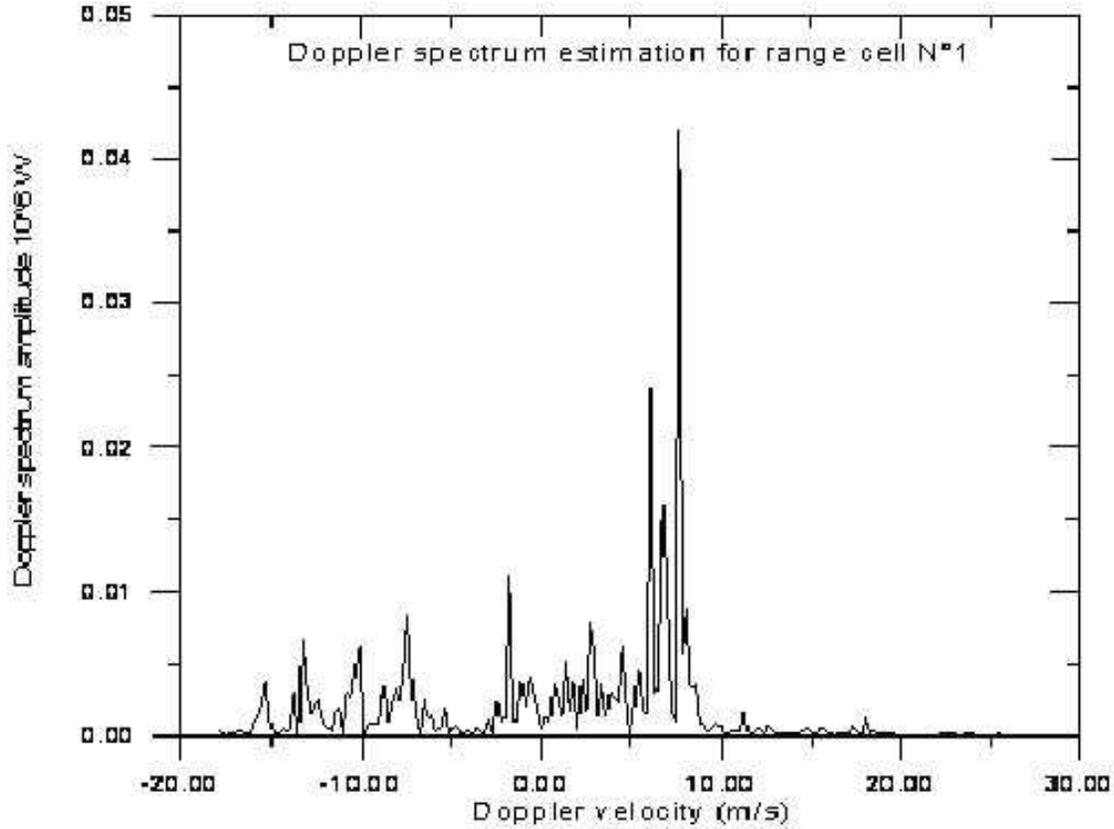}}
\caption{The received Doppler spectrum}
\end{figure}

\textbf{}

The received Doppler spectrum represents the power spectrum density of the received signal for a detection volume. Furthermore, without taking into account the noise power, the total echo power is given by the zeroth moment function, [1], [2]:

\begin{equation} \
P=\int _{}^{}\, S(v)\, dv
\end{equation}

The mean wind velocity where the first normalized moment is [2]:

\begin{equation} \
\bar{v}=\frac{1}{P} \int vS(v)dv
\end{equation}

The spectrum width of the Doppler spectrum mean velocity or shear is given by the square root of the normalised second central moment [2]:

\begin{equation} \
\sigma _{v}^{2} =\frac{1}{P_{} } \int _{}^{}(v-\bar{v})^{2} {\kern 1pt} \, S(v)\, dv
\end{equation}

The Doppler spectrum \textit{Sz(f)} can be scaled according to the Doppler velocity as \textit{S(v)}, using the relationship between the velocity and the Doppler frequency shift where \textit{$\lambda $} is the wavelength of the emitted signal: [3]

\begin{equation} \
v=\left(\frac{\lambda }{2} \right)f
\end{equation}

In the same way, we may draw the relationship between the velocity width $w$ and the standard Doppler spectrum deviation $\sigma _{f} $as [1], [2]:

\begin{equation} \
w=\left(\frac{\lambda }{2} \right)\sigma _{f}
\end{equation}

\subsection{Time domain estimation}

This approach is based on the estimation of the complex autocorrelation function of the radar signal [1]. A complex and stationary random process representing the time series of the radar signal (\textit{Z}(\textit{I},\textit{Q})), (I: in phase and Q: quadratic components) sampled at the pulse repetition time \textit{Ts} can be written as:

\begin{equation} \
Z(kT_{S} )=I(kT_{S} )+jQ(kT_{S} )
\end{equation}

Considering that these signals are statistically independent, thus the autocorrelation function is as follows [1]:

\begin{equation} \
R_{ZZ} (T_{s} )=\frac{1}{m} \sum _{k=0}^{m-1}Z^{*} (kT_{s} )Z((k+1)T_{s} )
\end{equation}

Where \textit{m} is the number of pulses considered. Assuming prior estimation of the noise power \textit{N} the total power can be estimated b [1], [6], [7]:

\begin{equation} \
\hat{P}=\frac{1}{m} \sum _{k=1}^{m}\left|Z(kT_{S} )\right|^{2} -N_{}
\end{equation}

The estimate of the mean wind velocity and its variance is given respectively by the estimators [1], [3], [6] and [7]:

\begin{equation} \
\hat{v}_{PP} =\frac{\lambda }{4\pi T_{S} } \arg \left[R_{ZZ} (T_{S} )\right]
\end{equation}

and

\begin{equation} \
\hat{w}_{PP}^{2} =\frac{\lambda ^{2} }{8\pi ^{2} T_{S}^{2} } \left[1-\frac{R_{ZZ} (T_{S} )}{\hat{P}} \right]
\end{equation}

The spectrum width of the mean wind velocity is obtained straightforward by the square root of the variance. The subscript (PP) stands for pulse-pair method.

\subsection{Frequency domain estimation}

Another technique for estimating the mean velocity$\hat{v}$, the variance$\sigma _{v}^{2} $ and the received Doppler spectrum width,$\hat{w}$ consists of estimating the power spectrum density via the Discrete Fourier Transform [5], [7], [11] as follows:

\begin{equation} \
\hat{V}_{FT} =\frac{\lambda }{2\hat{P}T_{S} } \sum _{k=-\frac{M}{2} }^{\frac{M}{2} -1}S_{Z} (k).\left(\frac{k}{M-1} \right)
\end{equation}

and:

\begin{equation} \
\hat{W}_{FT}^{2} =\frac{\lambda ^{2} }{4\hat{P}T_{S}^{2} } \sum _{k=-\frac{M}{2} }^{\frac{M}{2} -1}S_{Z} (k).\left(\frac{k}{M-1} +2\frac{\hat{V}_{FT} T_{S} }{\lambda } \right) ^{2}
\end{equation}

With,$S_{Z} (k)$is the power spectrum density.

These estimated values which are indexed (FT), are referred to as Fourier method.

This spectrum estimation is normalized by the total mean power,\textit{ P}, dealt with as a probability all over the considered bandwidth [4].

The advantage of the time implementation domain for the PP estimation is that it is less time-consuming than all the methods which require the use of the Discrete Fourier Transform (DFT) [6], [7].

The radar-reflected weather spectrum will not probably contradict these assumptions, but the presence of the inverted clutter mode (change of direction) can biase the mean spectrum estimation [6].

The DFT could be calculated with a fast DFT algorithm (FFT). However, the DFT has two disadvantages, which are inherent to its approach. The first one is that the frequency resolution is limited by the reverse of the width of the recorded samples. The second one implies the use of length-limited samples for the representation of infinite range signals. If we consider that the sequences will be worthless out of the finished interval, it means that a data windowing will be imposed. This windowing is equivalent to the data multiplication by a rectangular window of unit amplitude. As far as the frequency domain is concerned, the result is similar to a convolution of a spectrum with a \textit{sinc} function. [6][7].

This phenomenon is known as `spectrum losses' because the signal energy is not represented in the whole frequency domain [12][13].

These limitations of the DFT will arouse problems with short time sequences that are used in airborne Doppler radars [7]. As far as our study is concerned, it is not the case since we have adopted a ground-based pulse Doppler radar for our calculations.

\section{Ramanujan Fourier transform}

A In this section, notions of the Ramanujan\textit{ }sums, \textit{cq}(\textit{n}), have been dealt with. These sums may be defined as the \textit{nth} power of the \textit{qth }primitive roots of the unit [9].

\begin{equation} \
\begin{array}{l} {c_{q} (n)=\sum _{p=1}^{q}\exp (2i\pi \frac{p}{q} n) } \\ {\, \, \, \, \, \, \, \, \, \, \, (p,q)=1} \end{array}
\end{equation}

Where (\textit{p}, \textit{q}) =1 means that p and q are coprimes

And the coefficients, \textit{cq}(\textit{n}),  are sums from a set of \textit{ep}(\textit{n}) characters

With:

\begin{equation} \
e_{p} (n)=\exp (2i\pi \frac{p}{q} n)
\end{equation}

That the sums introduced by Ramanujan will play a fundamental role in the projection of arithmetic sequences,\textit{ x}(\textit{n}).

\begin{equation} \
x(n)=\sum _{q=1}^{\infty }x_{q} c_{q} (n)
\end{equation}

Glancing at the above equation, one can easily notice that infinite series having $q\to \infty $ may lead us to Fourier series [9], [10]. Moreover, the Discrete Fourier Transform takes a finite \textit{q} value. The arithmetic function,$\sigma (n)$, sum of \textit{n }dividers, can be written with RFT coefficients as $\sigma _{q} =\frac{\pi ^{2} n}{6} \frac{1}{q^{2} } $, consequently

\begin{equation} \
\sigma (n)=\frac{\pi ^{2} n}{6} \left\{1+\frac{(-1)^{n} }{2^{2} } +\frac{2\cos (2n\pi /3)}{3^{2} } +\frac{2\cos (n\pi /2)}{4^{2} } +.....\right\}
\end{equation}

For mean value functions, \textit{x}(\textit{n}), we have :

\noindent

\begin{equation} \
\begin{array}{l} {A_{v} (x)=\lim \frac{1}{t} \sum _{n=1}^{t}x(n) } \\ {\, \, \, \, \, \, \, \, \, \, \, \, \, \, \, \, \, \, t\to \infty } \end{array}
\end{equation}

Getting the inversion formula [9]

\begin{equation} \
x_{q} =\frac{1}{\phi (q)} A_{v} (x(n)c_{q} (n))
\end{equation}

From [9], the formula can be generalized. In what follows, the indexed coefficient,\textit{ xq,,} represent the Ramanujan Fourier Transform (RFT). The result is a multiplicative property of the Ramanujan~sums (coefficients):

\begin{equation} \
\noindent c_{qq'} (n)=c_{q} (n)c_{q'} (n)$ if $(q,q')=1
\end{equation}

Where the orthogonality property is:

\begin{equation} \
\sum _{n=1}^{q}c_{q}^{2} (n)=q\phi (q)
\end{equation}
The Ramanujan coefficients have been evaluated with the number theory functions. (\textit{q}, \textit{n}) represent the greater common divider of  \textit{q} and \textit{n.} Factoring a number into prime numbers, one may transcript q and \textit{n} as:

\begin{equation} \
\noindent q=\prod _{i}q_{i}^{\alpha _{i} }  (q_{i} prime)                                           \
\end{equation}

\begin{equation} \
\noindent n=\prod _{k}n_{k}^{\beta _{k} }  (n_{k} prime)
\end{equation}

We can also find the number $\phi (q)$of the irreducible fraction of denominator,\textit{ q}, also called Euler Totient function

\begin{equation} \
\phi (q)=q\prod _{i}(1-\frac{1}{q_{i} } )
\end{equation}

And Moebius function$\mu (n)$, with which these prime numbers could be obtained, is defined as

\begin{equation} \
\mu (n)=\left\{\begin{array}{c} {0\, \, if\, \, n\, \, \, {\rm content\; squared\; value}\, \, \beta _{k} >1\, \, \, \, \, \, \, \, \, \, \, \, \, \, \, \, \, \, \, \, \, \, \, \, \, \, \, \, \, \, \, \, \, \, \, \, \, \, \, \, \, \, } \\ {1\, \, {\rm if}\, \, n=1,\, \, \, \, \, \, \, \, \, \, \, \, \, \, \, \, \, \, \, \, \, \, \, \, \, \, \, \, \, \, \, \, \, \, \, \, \, \, \, \, \, \, \, \, \, \, \, \, \, \, \, \, \, \, \, \, \, \, \, \, \, \, \, \, \, \, \, \, \, \, \, \, \, \, \, \, \, \, \, \, \, \, \, \, \, \, \, \, \, \, \, \, } \\ {{\rm (-1)}^{{\rm k}} \, \, {\rm if}\, \, n\, \, {\rm is\; the\; product\; of\; k\; prime\; number}\, \, \, \, \, \, \, \, \, \, \, \, \, \, \, \, \, \, \, \, \, \, \, \, \, \, \, } \end{array}\right.
\end{equation}

From [9], the Ramanujan sums have been evaluated as follows:

\begin{equation} \
c_{q} (n)=\mu \left(\frac{q}{(q,n)} \right)\frac{\phi (q)}{\phi \left(\frac{q}{(q,n)} \right)}
\end{equation}

It should be noted that the sequences, \textit{cq}(\textit{n),} are periodic.

Having to face serious problems when calculating the mean wind velocity and its spectrum width with the Ramanujan-Fourier method, we had to calculate the spectrum density of the power deriving from the Discrete Fourier Transform (DFT) [12], [13].

\begin{equation} \
\hat{f}=\sum _{i}f_{i}  S_{RFT} (f_{i} )/\sum _{i}S_{RFT} (f_{i}  )
\end{equation}
And

\begin{equation} \
\hat{w}^{2} =\sum _{i}(f_{i} -\hat{f})^{2}  S_{RFT} (f_{i} )/\sum _{i}S_{RFT} (f_{i} )
\end{equation}
With

\begin{equation} \
\hat{v}=\lambda /2\hat{f}
\end{equation}

With $\hat{v}$ and $\hat{w}$ as being the estimated mean velocity and its spectrum width respectively, and $S_{RFT} (f_{i} )$ the Ramanujan spectrum of the complex series,\textit{ Z}(\textit{I}, \textit{Q}).

\section{Results and comments}

The data used are taken by a pulse Doppler radar WSR-88D at the state of Memphis Tennessee in July of 1997; they contains samples data I, Q, Azimuth, Elevation, Prt (pulse repetition times), Time (UNIX time).

The real data used are filtered by using IIR elliptic filter with fourth poles order, in order to reject the clutters and noises. This phase is very important to get only the weather spectrum modes.

We give in the figure 2 and figure 3, the representation of the complex series I and Q and their Doppler spectrum, for the Range cell n${}^\circ$1.

\begin{figure}[ht]
\centerline{\includegraphics[width=12.0truecm,clip=]{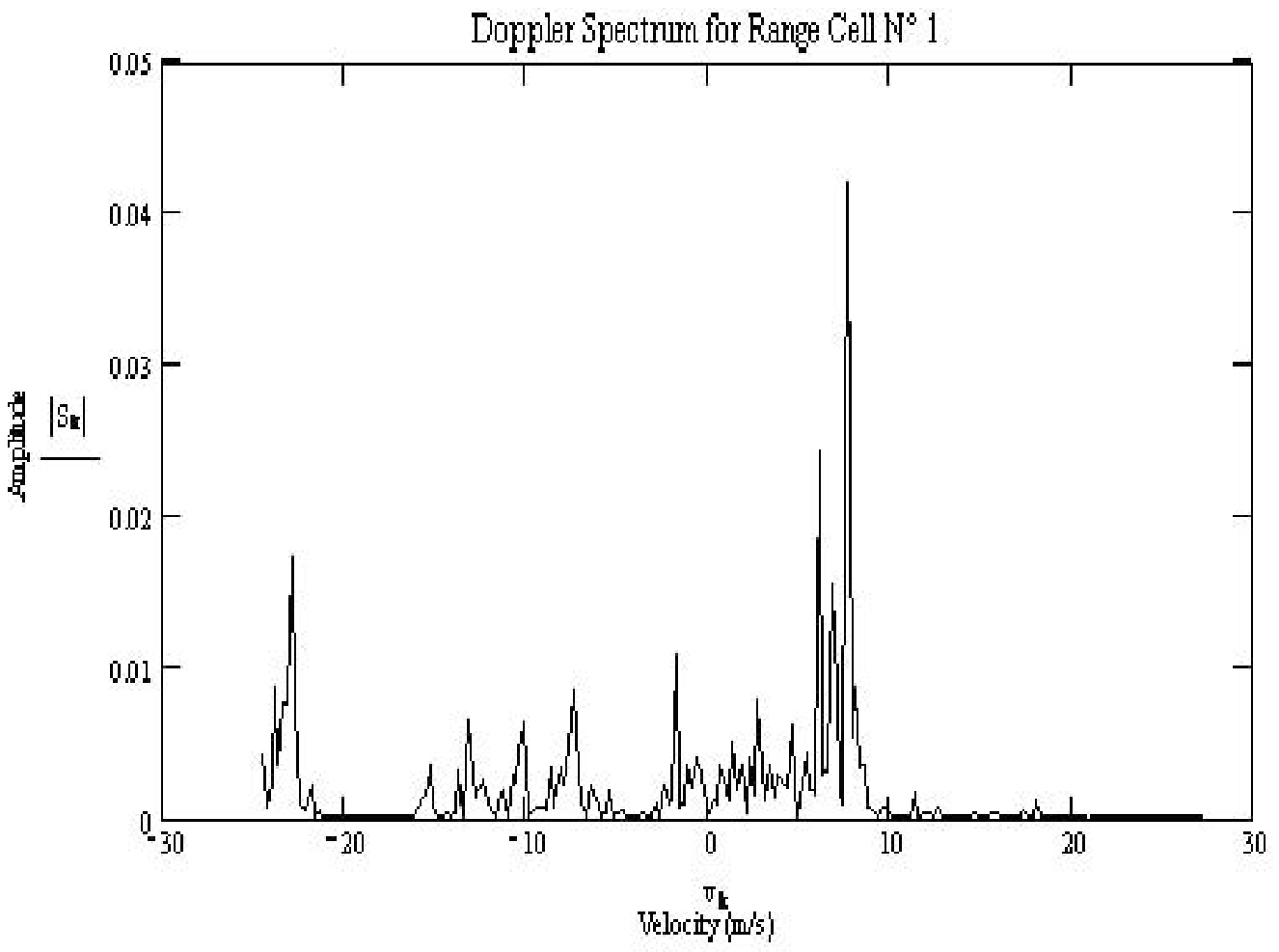}}
\caption{Time series I and Q}
\end{figure}

\begin{figure}[ht]
\centerline{\includegraphics[width=12.0truecm,clip=]{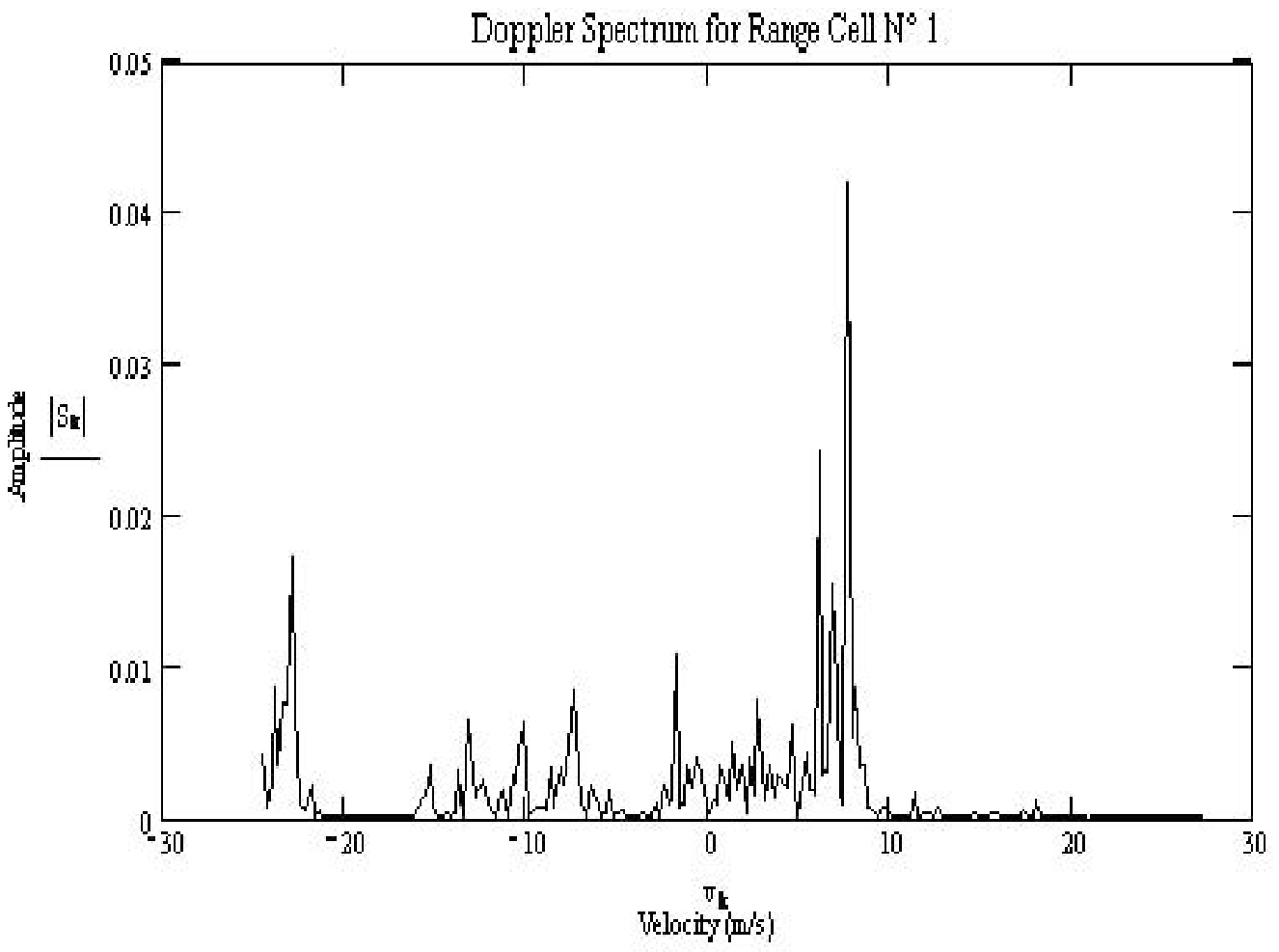}}
\caption{Doppler spectrum}
\end{figure}

\noindent

The three methods (pulse-pair, FFT and Ramanujan Fourier RFT), used for the estimation of the Doppler spectrum parameters of the wind perturbation led to the results represented in figures 5 and 6.

The algorithm based in the time domain (PP), is a simple method to program, because we have considered only the calculation of autocorrelation function of the received weather radar signals$Z(iT_{s} )$. It has less calculation time (13 ms) comparatively to RFT (16 ms), and much less to FFT (26ms), and converges at the first iteration (see table 1).

The estimation made by this method for the mean radial velocity of the wind Doppler spectrum, is very close to those provided in PPI form (Plan Position Indicator) by R.J. Keeler.

The variance $\sigma _{V}^{2} $ and the spectral width $w$ estimated by the pulse-pair method are weaker compared to those of the FFT and RFT methods.  The main disadvantage of this method is that the results related to the autocorrelation function are not easy to interpret, contrary to those of the spectral methods which are distributed over several frequencies (periodograms), [6], [13].

Moreover, the Ramanujan estimations which have been produced by the Ramanujan Fourier Transform (RFT), for the PSD calculation, are almost better to those provided by the FFT algorithm. The results are not different and the estimated values oscillate closely to those estimated about pulse-pair method. One advantage of using this method in calculations -- compared to FFT -- is that it is less time-computing, since it requires reduced samples only (see table 1).

\noindent

\noindent

\begin{figure}[ht]
\centerline{\includegraphics[width=15.0truecm,clip=]{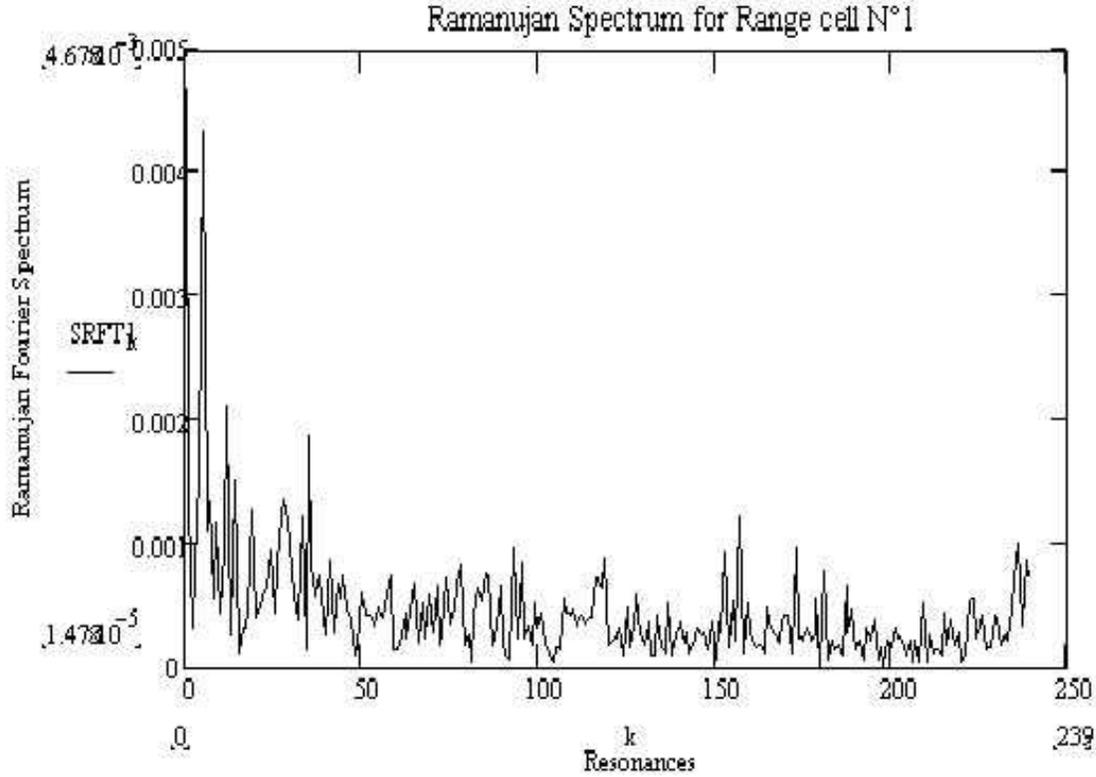}}
\caption{Ramanujan Spectrum for the range cell of N${}^\circ$1}
\end{figure}

\noindent

\noindent

\noindent
\begin{figure}[ht]
\centerline{\includegraphics[width=17.0truecm,clip=]{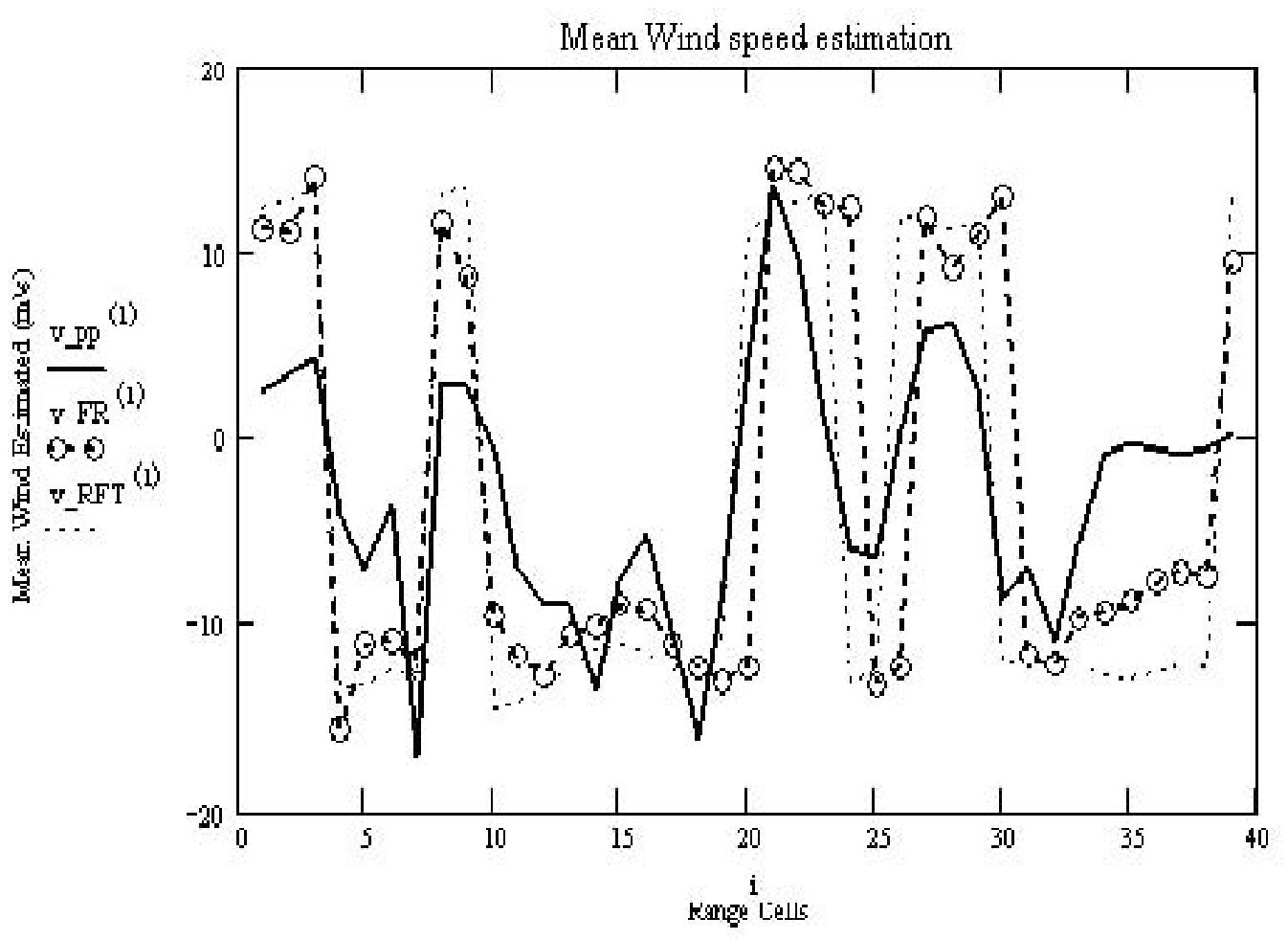}}
\caption{Mean wind velocity estimation}
\end{figure}

\noindent

In addition it is seen well (in the figure 5), these estimations done by Ramanujan Fourier algorithm are however better to those estimated by FFT method and non divergent. The disadvantage of this method lies in the execution time which is long to pulse-pair one. This long calculation time is due to the use of various coefficients for the \textit{PSD} estimation and which takes more calculation time. To reduce this execution time, there is a real time signal processing technique which uses electronic implementation circuit, Digital Signal Processor [11].  So, the estimation of the mean wind velocity and width done by the RFT algorithm will be improved and made in real time without any delay.

Consequently, the time - used for estimating the weather parameters (velocity, deviation/width) and forecasting different phenomena by the Ramanujan Fourier algorithm -- might be easily reduced.

\noindent

\noindent

\begin{figure}[ht]
\centerline{\includegraphics[width=17.0truecm,clip=]{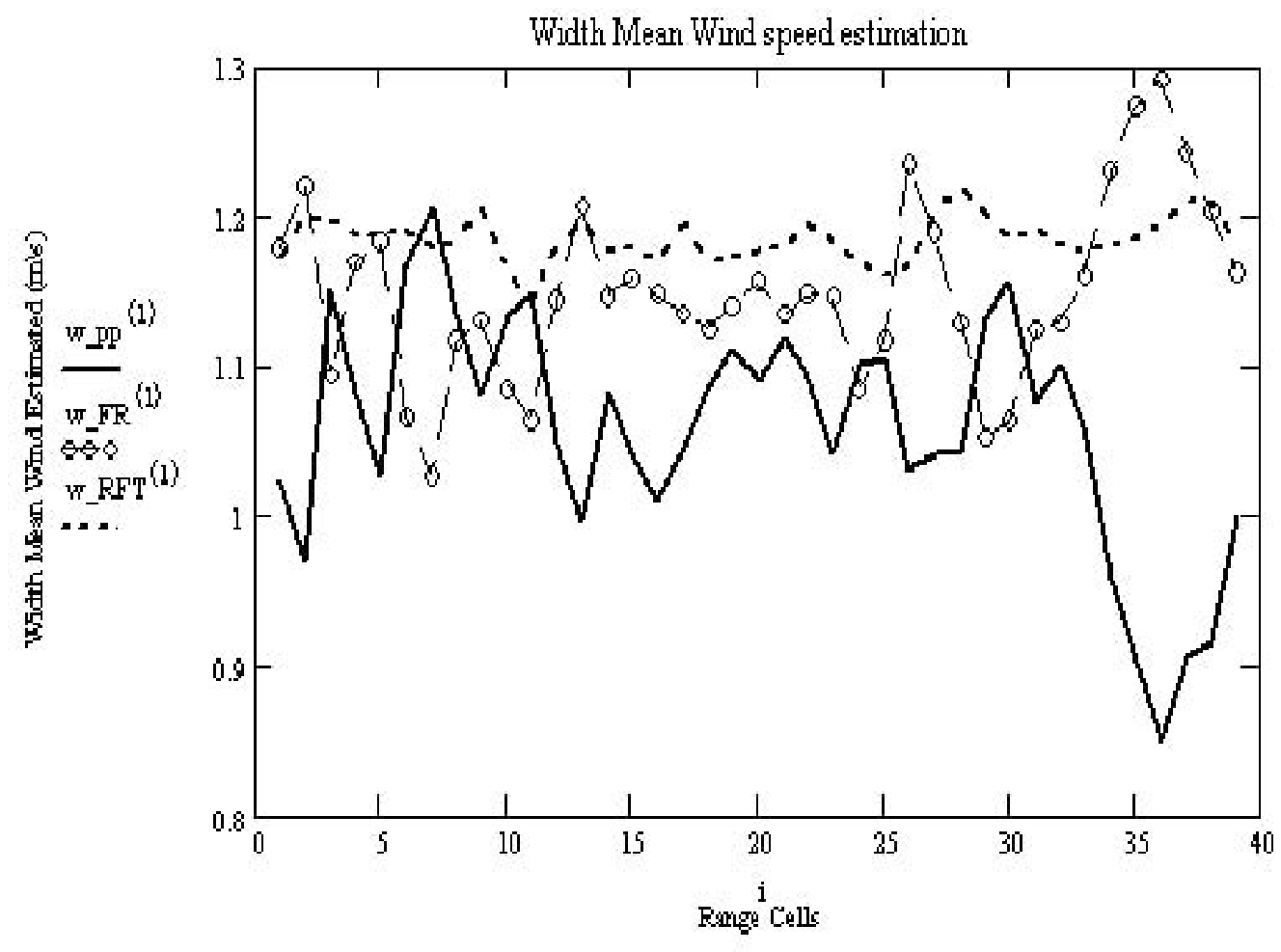}}
\caption{Spectral mean wind width estimation}
\end{figure}

\newpage

\noindent Table 1 : Calculation performances of the three algorithms
\vskip 5mm
\begin{tabular}{|p{1.2in}|p{1.in}|p{.8in}|p{.8in}|} \hline
Methods  & V[m/s] & W[m/s] & Calculation \\
&&& time (ms) \\ \hline
Pulse pair & 13.64 -- - 16.96 & 1.20  --  0.85 & 13 \\ \hline
FFT & 14.60 --   15.00 & 1.29  --  1.02 & 26 \\ \hline
Ramanujan Fourier & 13.65 --   14.40 & 1.22  --  1.14 & 16 \\ \hline
\end{tabular}
\vskip 5mm

It can be noted that the calculation performances of the three algorithms - which have been used to estimate the Doppler spectrum of a wind perturbation detected by pulsed Doppler ground-based radar -- have been mentioned in table 1.  The expected values, provided by the Ramanujan Fourier algorithm - for the mean wind and variance estimation - are not inferior to than those provided by the FFT and Pulse Pair ones. Meanwhile, it might be noted that the calculation time for a range cell is shorter than the one used by a classical FFT 16 ms Fourier algorithm (see table 1). The time-estimating algorithm has the lowest calculation time but its results are much more difficult to interpret because of the use of autocorrelation functions.

\noindent
\noindent

\section{Conclusion}

Using the Ramanujan Transform might be useful to identify a low magnitude spectrum. This transform is based on the prime numbers theory in relation with Moebius or Mangoldt functions.

This technique (RFT) permits faster calculations than the Discrete Fourier Transform since the calculations are made upon samples only; with respect to Moebius function. This technique takes only coprime resonances (p,q)=1. We might consider only a reduced number of samples because the estimation is made near zero frequency)

The spectrum estimation -- particularly after filtering I and Q samples (complex series Z) of the weather data -- seems to be very useful, knowing that the Ramanujan Transform facilitates the distinction of low-level weather modes like low intensity winds or rains.

Applying further notions of the number theory or even mathematics to the weather signal processing field may help obtain much more significant results and contribute for a better exploration of this wide and exciting field.

\noindent

\vskip 1cm

\noindent \textbf{Acknowledgments}

The authors would like to thank Dr. Michel Planat from the Laboratoire de Physique et Métrologie CNRS - Université de Franche-Comté France UMR6471, for its considerable contribution to the elaboration of this work. The authors express their gratitude's to R.J. Keeler and S. Ellis from NCAR country-regionUSA, which have provided a real data taken from place City Memphis State Tennessee with (WSR-88D) operational weather radar on July of 1997.

\noindent

\noindent

\vskip 1cm

\noindent \textbf{References}

\begin{enumerate}
\item [{[1]}] RJ Serafin - Skolnik, M., Meteorological Radar, Radar Handbook, 2nd Ed,1990.

\item  [{[2]}] R. J. Keeler and R. E. Passarelli, Signal Processing for Atmospheric Radars, (A90-39376 17-47),Radar in Meteorology, Boston, MA, American Meteorological Society, 1990, 199-229.

\item  [{[3]}] D. S. Zrni and J.T. Lee, Pulsed Doppler Radar Detects Weather Hazards to Aviation, Journal of Aircraft, Vol. 19, No 2, 1982, 183-190.

\item  [{[4]}] R. D. Palmer, Signal Processing Project, Department of Electrical Engineering University of Nebraska-Lincoln, spring 2002.

\noindent \item  [{[5]}] E. S. Chornoboy, Optimal mean velocity estimation for Doppler weather radars, IEEE transactions on geoscience and remote sensing, vol. 31, 1993, 575-586.

\noindent \item  [{[6]}] M. Lagha and M. Bensebti, "Performances Comparison of Pulse-Pair and 2-Step Prediction Algorithms for the Doppler Spectrum Estimation," Multidimensional Systems and Signal Processing, DOI: 10.1007/s11045-007-0030-3, 2008.

\noindent \item  [{[7]}] M. Lagha, H. M. Benteftifa and S. Boukraa, Contribution to the spectral estimate of the mean speed and variance of a windshear at low altitude- 1st International Symposium on Electromagnetism, Satellites, 2005.

\noindent \item  [{[8]}] S S. Abeysekera, Performance of Pulse-Pair Method of Doppler Estimation, IEEE Trans on Aerospace and Electronic Systems, 34 : 2, 1998, 520-531.

\noindent \item  [{[9]}] M.Planat, H. Rosu, and S. Perrine, Ramanujan sums for signal processing of low-frequency noise, Phys. Rev. E 66, 2002.

\noindent \item  [{[10]}] S Samadi, M.O Ahmad, and M. N. S. Swamy, Ramanujan Sums and Discrete Fourier Transforms, IEEE Signal Processing Letters, vol.12, no. 4, April 2005, 293-296.

\noindent \item  [{[11]}] Monakov, A.A.; Blagoveshchensky, D.V., A Method of Spectral Moment Estimation, IEEE Ttrans. Geosci and Remote Sens, vol. 37, no. 2, March 1999, 805-810.

\noindent \item  [{[12]}] V.M. Melnikov, D.S. Zrniç, Estimates of Large Spectrum Width from Autocovariances, Journal of Atmospheric and Oceanic Technology, 2004, 969-974.

\noindent \item  [{[13]}] J. Lee, Doppler Moment Estimation in a Weather Radar, International Journal of Electronics, Taylor \& Francis, 2002, 583 - 592.

\end{enumerate}

\end{document}